\newcommand{\lab}[1]{\label{#1}}                % hides labels

\documentclass[11pt]{article}
\usepackage{amsmath,amssymb,mathrsfs}
\usepackage{amssymb}
\usepackage[left=2.5cm,top=2.5cm,right=2.5cm,bottom=2.5cm]{geometry}
\usepackage{epsfig}
\usepackage[usenames,dvipsnames]{color}

\begin{document}
\newtheorem{theorem}{Theorem}[section]
\newtheorem{lemma}[theorem]{Lemma}
\newtheorem{definition}[theorem]{Definition}
\newtheorem{conjecture}[theorem]{Conjecture}
\newtheorem{proposition}[theorem]{Proposition}
\newtheorem{algorithm}[theorem]{Algorithm}
\newtheorem{corollary}[theorem]{Corollary}
\newtheorem{observation}[theorem]{Observation}
\newtheorem{problem}[theorem]{Open Problem}
\newcommand{\noin}{\noindent}
\newcommand{\qed}{\ \hfill \rule{1ex}{1ex}} %rightjustified box
\newcommand{\ind}{\indent}
\newcommand{\om}{\omega}
\newcommand{\N}{{\mathbb N}}
\newcommand{\R}{{\mathbb R}}
\newcommand{\Z}{{\mathbb Z}}
\newcommand{\RR}{{\mathscr{R}}}
\newcommand{\E}{\mathbb E}
\newcommand{\Prob}{\mathbb{P}}
\newenvironment{proof}{{\noin \bf Proof}: }{\qed}

\newcommand{\bee}{\begin{equation}}
\newcommand{\ee}{\end{equation}}
\newcommand{\bea}{\begin{eqnarray}}
\newcommand{\eea}{\end{eqnarray}}
\newcommand{\non}{\nonumber}
\newcommand{\bean}{\begin{eqnarray*}}
\newcommand{\eean}{\end{eqnarray*}}
\newcommand\eqn[1]{(\ref{#1})}
\newcommand{\bel}[1]{\bee\lab{#1}}

\newcommand{\hatk}{{\bf k}}
\newcommand{\haty}{{\bf y}}
\newcommand{\hatz}{{\bf z}}
\newcommand{\hatw}{{\bf w}}
\newcommand{\tildez}{\tilde{\bf z}}

\title{Almost all 5-regular graphs have a 3-flow}

\author{
Pawe{\l} Pra{\l}at\thanks{Supported by NSERC.}\\ {\small Department of Mathematics}\\ {\small Ryerson University}\\ {\small Toronto  ON M5B 2K3}\\
{\small  Canada.}
\and Nick Wormald\thanks{Supported by an ARC Australian Laureate Fellowship.}\\
{\small School of Mathematical Sciences}\\{\small Monash University VIC 3800}\\
{\small Australia.} 
}
\date{}

\maketitle

\begin{abstract}
Tutte conjectured in 1972  that every 4-edge connected graph has a nowhere-zero 3-flow. This has long been known to be equivalent to the conjecture that every 5-regular 4-edge-connected graph has an edge orientation in which every out-degree is either 1 or 4. We show that the assertion of the conjecture holds asymptotically almost surely for random 5-regular graphs. It follows that the conjecture holds for almost all 4-edge connected 5-regular graphs.
\end{abstract}

\section{Introduction}

A {\em nowhere-zero $3$-flow} (sometimes simply called a 3-flow) in an undirected graph $G=(V,E)$ is an orientation of its edge set $E$ together with a function $f$ assigning a number $f(e) \in \{1,2\}$ to every $e\in E$ such that the following is satisfied. For every vertex $v \in V$,
$$
\sum_{e \in D^+(v)} f(e) - \sum_{e \in D^-(v)} f(e)=0,
$$
where $D^+(v)$ is the set of  all edges oriented away from   $v$,  and $D^-(v)$ is the set of all edges oriented towards $v$.

A well known conjecture of Tutte from 1972 (see e.g.\ Bondy and Murty~\cite{BM} (Open Problem 48) and Jensen and Toft~\cite[Section 13.3]{JT})  asserts that every $4$-edge-connected graph  admits a nowhere-zero $3$-flow. This conjecture is still open. For long, it was  not even known whether or not there is a fixed $k$ such that every $k$-edge connected graph has a nowhere-zero $3$-flow (known as the weak 3-flow conjecture of Jaeger).  
Weaker versions, for $k\ge c \log_2 n$ for $n$-vertex graphs, proved by Alon, Linial and Meshulam~\cite{ALM} and Lai and Zhang~\cite{LZ}. Recently, the weak 3-flow conjecture was settled by Thomassen~\cite{Thomassen}, who proved that every $8$-edge-connected graph admits a nowhere-zero 3-flow. This was subsequently improved to $k=6$ by Lov\'{a}sz, Thomassen, Wu, and Zhang~\cite{LTWZ}.

It is known (see, e.g., Seymour~\cite{Se}) that a graph admits a nowhere-zero $3$-flow if and only if it has a nowhere-zero flow over $\Z_3$, or equivalently,  an edge orientation in which the difference between the out-degree and the in-degree of every vertex is divisible by $3$. It has also long been known (see~\cite{BM} and~\cite{JT}) that it is enough to prove the conjecture for $5$-regular graphs. Thus, Tutte's conjecture has the following equivalent form.

\begin{conjecture}[Tutte]\label{c11}
Every $4$-edge-connected $5$-regular graph has an edge orientation in which every out-degree is either $1$ or $4$.
\end{conjecture}

In this paper, we show that Tutte's conjecture holds for almost all 5-regular graphs. To state this precisely, we say that a property of a probability space indexed by $n$ holds {\em a.a.s.} (asymptotically almost surely)  if the probability that it holds tends to $1$ as  $n$ tends to $\infty$ (with $n$ restricted to being even for odd-degree regular graphs). Using the small subgraph conditioning method of Robinson and Wormald~\cite{RW} (see~\cite{Wo}) we show the following.

\begin{theorem}\lab{thm:main}
A random $5$-regular graph $G_n$ on $n$ vertices a.a.s.\ admits  a nowhere-zero 3-flow, that is, an edge orientation in which every out-degree is either $1$ or $4$.
\end{theorem}

Since it is well known that almost all 5-regular graphs are 5-edge-connected (see e.g.~\cite{Wo}), it follows that almost all 4-edge-connected $5$-regular graphs have a nowhere-zero 3-flow.

Jaeger~\cite{Ja} generalised Conjecture~\ref{c11} by conjecturing that for any integer $p \geq 1$, the edges of every $4p$ edge-connected graph can be oriented so that the difference between the out-degree and the in-degree of every vertex is divisible by $2p+1$. Similar to Tutte's conjecture, it is known that the general case can be reduced to the $(4p+1)$-regular case. Alon and Pra{\l}at~\cite{AlonPralat} showed that the assertion of Jaeger's conjecture holds for almost all $(4p+1)$-regular graphs, provided that $p$ is large enough.  (The lower bound for $p$ was not optimized, but it could not be reduced to $p=1$.)  The proof used methods quite different from the present paper,  involving an application of the Expander Mixing Lemma to an equivalent version of the conjecture.

\section{Proof of Theorem~\ref{thm:main}}

The pairing model for investigating properties of random regular graphs was instigated by Bollob{\'a}s~\cite{Bo}. This consists  of $dn$ points that are arranged in $n$ groups (called vertices) of $d$ each, arranged in pairs uniformly at random. The pairs induce a multigraph in the obvious way, and we refer to pairs as edges. This pairing model, called $\mathcal{P}_{n,d}$, is useful because  simple graphs occur with equal probabilities, and the probability that it is simple for fixed $d$ is bounded away from 0. Hence, to show that the random regular graph has a property a.a.s., it is enough to show that the random member of the multigraph corresponding to $\mathcal{P}_{n,d}$ a.a.s.\ has the same property or is non-simple. (See~\cite{Wo} for more information on this and other claims we make about $\mathcal{P}_{n,d}$.)  We will work with   orientations of (the pairs of) a pairing in $\mathcal{P}_{n,5}$ in which each vertex has in-degree 1 or 4. We call such orientations {\em valid}. Given   an orientation, vertices of in-degree 1 will be called \emph{in-vertices}, and those of out-degree 1 {\em out-vertices}, and each point contained in an edge  oriented towards an in-vertex, or away  from an out-vertex, is called \emph{special}. Moreover, a point is an {\em in-point} if the   edge containing it is oriented towards it, and an {\em out-point} otherwise.

Let $Y=Y(n)$ be the number of valid orientations of a random element of $\mathcal{P}_{n,5}$. It is easy to see that
$$
\E Y  = \frac { {n \choose n/2} 5^n (5n/2)!}{M(5n)},
$$
where
$$
M(s) = \frac {s!}{(s/2)! 2^{s/2}}
$$
is the number of perfect matchings of $s$ points. Indeed, there are $n \choose n/2$ ways to select in-vertices (since exactly half of the vertices must be such), $5^n$ ways to select one special point  in each vertex, which determines each point to be either in or out, $(5n/2)!$ ways to pair up the points so that each ``in'' is paired with an ``out'', and $M(5n)$ pairings in total. Using Stirling's formula  $s! \sim \sqrt{2 \pi s} (s/e)^s$, we get 
\bel{firstmoment}
\E  Y  = \frac {n! 5^n (5n/2)!^2 2^{5n/2}}{(n/2)!^2 (5n)!} \sim \left( \frac {25}{8} \right)^{n/2} \sqrt{5}.
\ee

This tells us that there are plenty of valid orientations per pairing, on average. To show that pairings a.a.s.\  have at least one valid orientation, i.e.\ that $\Prob(Y>0) \sim 1$, a common method would be to estimate $\E  Y(Y-1) $, show that it is asymptotic to $(\E Y)^2$, and then apply Chebyshev's inequality. As we shall see, this fails in the present case, but only just,  as there is a constant factor discrepancy in the asymptotics. Under such circumstances, we can hope to apply the small subgraph conditioning method~\cite{Wo}. Here the first step is again to compute $\E  Y(Y-1)$, and then also some joint moments of $Y$ with short cycle counts (and then hope for the best).
 
To estimate $\E  Y(Y-1)$, consider any two orientations of the same 5-regular graph. Suppose that precisely  $k$ vertices are in-vertices in both orientations, and that, of these, precisely $k_{11}$   have the same special   point in both orientations. Since the first orientation induces $n/2$  in-vertices, exactly $n/2-k$ vertices   are in-vertices in the first orientation but out-vertices in the second one.  Of these,  suppose that $k_{10}$  ($k_{10} \le n/2-k$) have the   two special points coinciding. Similarly, there are $k$ vertices that are out-vertices in both orientations; suppose that $k_{00}$ of them  ($  k_{00} \le k$) have the two special points coinciding. Finally, there are $n/2-k$ vertices that are out- in the first and in- in the second orientation; suppose that $k_{01}$ of them  ($ k_{01} \le n/2-k$) have coinciding special points.

It turns out that there are no additional restrictions on these parameters, other than integrality and non-negativity. We define
$$
I = I(n)= \left\{ (k, k_{00}, k_{01}, k_{10}, k_{11}) \in \N_0^5 :  k \le \frac {n}{2}, \  \max\{ k_{00},k_{11}\} \le k, \ \max\{ k_{01},k_{10}\} \le \frac {n}{2} - k  \right\} 
$$
where $\N_0 = \N\cup \{0\}$.
Fix $\hatk = (k, k_{00}, k_{01}, k_{10}, k_{11}) \in I$. We next calculate the number of configurations, i.e.\   pairings with two given orientations, corresponding to this vector. There are
$$
\frac {n!}{k_{00}! k_{01}! k_{10}! k_{11}! (k-k_{00})! (k-k_{11})! (n/2-k-k_{01})! (n/2-k-k_{10})!}
$$
ways to partition the vertices into the eight groups. There are then
$$
5^{k_{00}+k_{01}+k_{10}+k_{11}} (5 \cdot 4)^{(k-k_{00})+(k-k_{11})+(n/2-k-k_{01})+(n/2-k-k_{10})} = 5^n \cdot 4^{n-k_{00}-k_{01}-k_{10}-k_{11}}
$$
ways to assign special points in the two orientations. Next, we need to pair (in,in)-points with (out,out)-points (where the first  ``in''  refers to the first orientation, and so on), and (in,out)- with (out,in)-points. The number of (in,in)-points   is equal to
$$
k_{11}+4k_{00}+3(k-k_{00})+(n/2-k-k_{01})+(n/2-k-k_{10}) = n+k+k_{00}+k_{11}-k_{01}-k_{10},
$$
and the same applies for (out,out). These two sets must be paired with each other.
Half of the  remaining points, or
$$
5n/2 - (n+k+k_{00}+k_{11}-k_{01}-k_{10}),
$$
will be in-out, and an equal number will be out-in.
Hence, there are 
$$
(n+k+k_{00}+k_{11}-k_{01}-k_{10})! (3n/2 - (k+k_{00}+k_{11}-k_{01}-k_{10}))! 
$$
ways to legally pair the points. The number of configurations is the product of the above factorials. To obtain the expected number of pairs of orientations, we must divide by the number $M(5n)$ of pairings. 
Putting $\hatz =\hatz(\hatk) = \hatk/n$ and applying Stirling's formula again, we can write
\bel{secondmoment}
\E Y(Y-1) =\sum_{\hatk \in I} r(\hatz)g(\hatz) \exp \big( nf(\hatz)\big),
\ee
where the various factors are defined as follows. 
The function $r$, 
which is the error factor in the applications of Stirling's formula, has the property that
$r=O(1)$ for all $\hatz$,  and $r\sim 1$ if $\hatz$ is bounded away from the boundary of 
$$
J:=\left\{ (z, z_{00}, z_{01}, z_{10}, z_{11}) \in \R_0^5 :  z \le \frac {1}{2},\ \max\{z_{00},z_{11}\} \le z,\ \max\{z_{01},z_{10}\} \le \frac {1}{2} - z  \right\}
$$
where $\R_0$ is the set of non-negative reals. With
$b=z+1+z_{00}-z_{01}-z_{10}+z_{11}$ and $h(x)=x \log  x$,  
$$
g=  
\frac{1}{ \sqrt{ 32} (\pi n)^{5/2}}\left(\frac{b(5-2b) }{ z_{00} z_{01} z_{10} z_{11} (z-z_{00}) (z-z_{11}) (1-2z-2z_{10}) (1-2z-2z_{01})}\right)^{1/2}
$$
from the polynomial factors in Stirling's formula, and 
\bean
f&= &\Big(  9/4-z_{00}-z_{01}-z_{10}-z_{11}\Big)\log 4+\log 5 -h(5)+h(5/2)  + h(b)+h(5/2-b)-h(z_{00})\\
&&-h(z_{01})-h(z_{10})-h(z_{11})-h(z-z_{00})-h(z-z_{11})-h(1/2-z-z_{01})-h(1/2-z-z_{10})
\eean
from the rest.
 
Note that we can extend the definition of $f$ continuously to the boundary of $J$ by defining $x\log x=0$ at $x=0$. Then $f$ achieves its maximum on $J$. Our next goal is to show that  $\tildez = (1/4,1/20,1/20,1/20,1/20)$ is the unique global maximum point of $f$ on $J$, since we can easily argue then that points far away from $\tildez$ give negligible contribution to~\eqn{secondmoment}. 

We first investigate stationary points in the interior of $J$. An algebraic manipulation package, such as the Maple we used, makes this easy. We find 
\bel{dz00}
\frac{\partial f}{\partial z_{00}}
= \log\frac{   (z-z_{00})(z+1+z_{00}-z_{01}-z_{10}+z_{11})}{2z_{00}(3-2z-2z_{00}+2z_{01}+2z_{10}-2z_{11})}.
\ee
Setting this equal to 0 gives $P_{00}=0$ where
$$
P_{00}=(z-z_{00})(z+1+z_{00}-z_{01}-z_{10}+z_{11})-2z_{00}(3-2z-2z_{00}+2z_{01}+2z_{10}-2z_{11}).
$$
Defining $P_{01}$ etc.\ similarly, and $P$ from $\partial f/\partial z $, we obtain five polynomials such that any local maximum must be a common zero of all five polynomials. Write $\RR(X,Y,x)$ for the resultant of two polynomials $X$ and $Y$ with respect to $x$. When $X=Y=0$, it is necessary that $\RR(X,Y,x)$=0. We find that
$$
P_6:= \RR(P_{00},P_{01},z_{10}) = 5z-5z_{00}-10z^2+10zz_{00}-10z_{01}z-150z_{01}z_{00},
$$
$$ 
P_6- \RR(P_{00},P_{10},z_{01}) =  10(-z_{10}+z_{01}))(z+15z_{00}).
$$
On the interior of $J$, we have $z>0$, so we may assume that $z_{01} = z_{10}$. Also,
$$
\RR\big(P_6,\RR(P_{01},P_{11},z_{10}),z_{01}\big) = -800z(z_{00}-z_{11})(-1+2z).
$$
On the interior $z < 1/2$,  so we may assume that $z_{00} = z_{11}$.

For any  polynomial $X$, let $X^*$ be the result of setting $z_{11}=z_{00}$ and $z_{10}=z_{01}$. We find
$$
P_7:=\RR(P^*_{00},P^*_{01},z_{01})= -60z^3-480z^2z_{00}-120zz_{00}-1500zz_{00}^2+2040z_{00}^2-1800z_{00}^3,
$$
$$
\RR\big(P_7,\RR(P^*,P^*_{00},z_{01}),z_{00} \big)=
  573308928\times 10^{11}\times z^8(4z-1)(13068z^2-6534z-109)( 1-2z)^2.
$$
Neither root of the quadratic factor is in $[0,1/2]$, so at any interior stationary point  we must have $z=1/4$.
Substituting  $z=1/4$ into $P_7=0$ gives
$$
 (20z_{00}-1) (96z_{00}^2-84z_{00}-1)=0.
$$
Again, the roots of the quadratic factor are out of range, and hence $z_{00}=1/20$ at any interior stationary point. Thus 
$z_{11}=1/20$ also. Substituting the known values into  $P^*_{01}$ gives 
$$
\frac{1}{20} (20z_{01}-1)(24z_{01}-23)
$$
and hence the unique   stationary point of $f$ in the interior of $J$ is $\tildez$.

We turn next to the boundary of $J$. First consider any point on the boundary at which $0<z<1/2$. Then $z_{01}+z_{10}\le 1-2z$ and $z_{00}+z_{11}\le 1 < 3/2-z$. Hence, for  $z_{00}$  tending towards 0 or $z$,   $\partial f/\partial z_{00}$ is dominated by the terms $-\log  z_{00}$ and $ + \log(z-z_{00})$, which tend to $+\infty$ and $-\infty$  respectively in the two cases.  This shows that for $0<z<1/2$, the  boundary points where $z_{00}=0$ or $z_{00}=z$ cannot contain a global maximum of $f$ on $J$. A similar observation applies also to show that $z_{11}$,  $z_{10}$ and $z_{01}$ cannot be at the extreme ends of their ranges. Hence, no such boundary point is a maximum of $f$.
 
 We are left with considering the boundary points where $z=0$ or $z=1/2$. 
 \medskip
 
\noindent
{\em Case  1. $z=0$.}
\smallskip

 Membership of $\hatz$ in $J$ then forces $z_{00}=z_{11}=0$. Substituting $z=z_{00}=z_{11}=0$ into $f$ gives a function $\bar f(z_{01},z_{10})$ with domain $[0,1/2]^2$. Setting its partial derivatives to 0 give two polynomials equal to 0. Their difference is
$$
 - 2(6z_{01}-11+6z_{10}) (-z_{10}+z_{01}).
$$
Since $z_{01}+z_{10}\le 1$, this is only 0 at $z_{01}=z_{10}$. We find $d \bar f(z_{10},z_{10})/dz_{10} = (3+4z_{01})^2/256z_{01}^2$, which is never 0. So no interior point can be maximum for $\bar f$.  This  leaves  the   boundary of its domain, where $z_{01}$ or $z_{10}$ is 0 or $1/2$.  If $z_{01}=0$,  then $\partial \bar f\partial z_{01}$  is large and positive, so there is no maximum there. By symmetry, it is a similar conclusion if $z_{10}=0$. On the other hand, along the boundary $z_{01}=1/2$, $d \bar f(1/2,z_{10})/dz_{10}=\log((2+z_{10})/4z_{10})$, which is always positive, and similarly for $z_{10}=1/2$. So the only possible local maximum on the boundary in this case is at $(0,0,1/2,1/2,0)$. Here $f=\log(5/8)$, whereas $f(\tildez)=\log(25/8)$.
 \medskip
 
\noindent
{\em Case  2. $z=1/2$.}
\smallskip

Membership of $\hatz$  in $J$ then forces $z_{01}=z_{10}=0$. It turns out that substituting $z=1/2$ and $z_{01}=z_{10}=0$ into $f$ produces exactly the function $\bar f (z_{00},z_{11})$, with the same domain. Hence Case~1 shows that the only local maximum of $f$ here is $\log(5/8)$.

In conclusion, $f$ has (at most) three local maxima, at $\tildez$, $(0,0,1/2,1/2,0)$ and $(0,1/2,0,0,1/2)$, and the first of these is the global maximum on $J$.

The next part of the proof consists of a routine computation and argument. Set $z=1/4+y$ and $z_{ij}=1/20+y_{ij}$ for each of the other four variables. Then expand $f(\hatz)$ about the global maximum point $\tildez$, to obtain  
\bel{fexp}
f=\log(25/8) +\haty^TB\haty +O(x^3) 
\ee
where $\haty=(y,y_{00},y_{01},y_{10},y_{11})$, $\haty^T$ denotes the transpose of $\haty$, 
$$
B=\frac{1}{10}
\left[
\begin{array}{r r r r r}
    -92&  33&-33&  -33&  33\\
     33&-117& -8&   -8&   8\\
    -33&  -8&-117&   8&  -8\\
    -33&  -8&   8&-117&  -8\\
     33&   8&  -8&  -8& -117
\end{array}
\right],
$$
and $x=x(\hatz) =||\haty||$, with $||\cdot||$ denoting  the $L^2$ norm (say). The error term in this expansion is valid  by Taylor's theorem provided that $x=o(1)$. 

A standard argument (with details given below) now allows us to estimate the summation in~\eqn{secondmoment}, and we obtain
\bel{second:mark2}
\E Y(Y-1) \sim \left(\frac{25}{8}\right)^n\frac{g(\tildez)(\pi n)^{5/2}}{\sqrt{|\det B|}}.
\ee
We find that $g(\tildez) = (5^5/2)(\pi n)^{-5/2}$ and $\det B = -3\cdot 5^6\cdot 7/4$, so~\eqn{firstmoment} gives
%\bel{second:final}
%\E Y(Y-1) \sim  \left(\frac{25}{8}\right)^n  \frac{25}{\sqrt{21}}.
%\ee 
$\E Y(Y-1) \sim \frac{25}{\sqrt{21}}  (25/8)^n.$
  Combining this with~\eqn{firstmoment}, we have 
\bel{ratio}
 \frac {\E \, Y^2 }{ (\E Y)^2 } \sim \frac {5}{\sqrt{21}}.
\ee

%One version of the above-mentioned standard argument is (briefly)  as follows.
We now justify~\eqn{second:mark2}. Let $J_0:= \{\hatz: x=o(n^{-2/5})\}$. For $\hatz\in J_0$, we have $r(\hatz)g(\hatz)\sim  g(\tildez)$ and $x^3=o(n^{-6/5})$. Thus
$$
 \sum_{\hatk: \hatk/n\in J_0} r(\hatz)g(\hatz) \exp \big( nf(\hatz)\big)
 \sim  \left(\frac{25}{8}\right)^ng(\tildez)\sum_{\hatk: \hatk/n\in J_0} e^{n\haty^TB\haty}.
 $$
The eigenvalues of $B$ are $(-37 \pm \sqrt{697})/4$, each with multiplicity 1, and $ -25/2$ with multiplicity 3. These are all less than $-2.6$. Hence $B$ is negative definite, and we have
$$
  \sum_{\hatk: \hatk/n\in J_0} e^{n\haty^TB\haty}\sim n^5\int_{J_0} e^{n\haty^TB\haty} d^n \haty.
$$
Here the factor $n^5$ arises from the change of variable from $\hatk$ to $\hatz=\hatk/n$, and the replacement of the summation by the integral can be justified by various elementary means. (One way is to rescale to $\hatw=\sqrt n \haty$ and observe that the summations for various $n$ amount to different Riemann sums for the integral of  the fixed function $ e^{ \hatw^TB\hatw}$. The regions of integration are different for different $n$ but the tails of the summation and the integral outside the regions are easily seen to be negligible.) A simple argument, for example   bounding the value of the integrand in terms of $||\haty||$, shows that the latter integral can be extended to all of $\R^5$ with no significant change, that is 
$$
\int_{J_0} e^{n\haty^TB\haty} d^n \haty \sim \int_{\R^5} e^{n\haty^TB\haty} d^n \haty =\frac{\pi^{5/2}}{ n^{5/2}\sqrt{|\det B|}}.
$$
Combining the last few conclusions gives
\bel{main}
 \sum_{\hatk: \hatk/n\in J_0} r(\hatz)g(\hatz) \exp \big( nf(\hatz)\big)
 \sim  \left(\frac{25}{8}\right)^n\frac{g(\tildez)(\pi n)^{5/2}}{\sqrt{|\det B|}}.
\ee

Recalling that $B$ is negative definite, and we conclude from~\eqn{fexp} that on the boundary of $J_0$, the value of $f$ is $f(\tildez)-\Omega(n^{-4/5})$. Since $f$ is a fixed function independent of $n$, and $\tildez$ is its global maximum, 
$$
\max_{\hatz\in J\setminus J_0}f(z)=  f(\tildez)-\Omega(n^{-4/5})  = \log(25/8) - \Omega(n^{-4/5}).
$$
Considering that $r$ and $g$ are polynomially bounded, it follows that each term in the summation in~\eqn{secondmoment} (where $f$ is multiplied by $n$) for which $\hatk/n\in J\setminus J_0$ is $(25/8)^ne^{-\Omega(n^{1/5})}$. There are a polynomial number of such terms, so their sum is likewise $(25/8)^ne^{ -\Omega(n^{1/5})}$. Thus the   terms in~\eqn{main} dominate the summation in~\eqn{secondmoment}, and we have~\eqn{second:mark2} as claimed.

The second main computation in the small subgraph conditioning method involves the random variable $X_k$ ($k \ge 1$) which is the number of cycles of length $k$ in $\mathcal{P}_{n,5}$.  It is known that for each $k \ge 1$, $X_1, X_2, \ldots, X_k$ are asymptotically independent Poisson random variables with 
\bel{lambda}
\E  X_k  \to \lambda_k := \frac {4^k}{2k}.
\ee
We are required to show, for each $k \ge 1$, that there is a constant $\mu_k$ such that
\bel{mudef}
 \frac {\E  (Y X_k) }{\E  Y} \to \mu_k 
\ee
and, more generally, such that the joint factorial moments satisfy
\bel{mugen}
\frac {\E (Y[ X_1]_{j_1} \cdots [X_k]_{j_k} ) }{\E  Y }  \to \prod_{i=1}^k \mu_i^{j_i}
\ee
for any fixed $j_1,\ldots, j_k$. (Here $[x]_k$ is the falling factorial.) We will derive a value of $\mu_k$ satisfying~\eqn{mudef} and will observe that essentially the same argument generalises easily to give~\eqn{mugen}. 

To evaluate $\E  (Y X_k)$, we  find the number of triples $(P,C,O)$ where $P$ is a pairing, $C$ a $k$-cycle of $P$ and $O$ an orientation of $P$, and then divide by $|\mathcal{P}_{n,5}|=  M(5n)$. 

The number of ways to choose the  pairs of (i.e.\ inducing the edges of) the cycle is $20^k[n]_k/2k$. The easiest way to see this is to choose the vertices of the cycle in one direction, starting at a canonical one ($[n]_k$), and then  an ordered pair of points   in each vertex ($5^k4^k$), and divide by the multiplicity of counting due to the canonical start and direction ($2k$).

We will count the triples $(P,C,O)$ which have $i$  vertices on $C$ with in-degree 2 {\em in $C$}.  These must be out-vertices in the orientation $O$. The number of vertices with out-degree 2 in $C$ must also equal $i$, and these must be in-vertices. The remaining $k-2i$ cycle vertices have  in-degree and out-degree  in $C$ both equal to 1.  These can be either in- or out-vertices. We now select the rest of the oriented pairing. There are $n - 2i \choose n/2-i$ ways to select the remaining in- and out-vertices,  and $3^{2i}$ ways to choose the special points of the vertices of $C$. (Note that the latter only needs to be done for vertices of in-degree 0 or 2 in $C$; vertices of in-degree 1 in the cycle have their special point already determined by the orientation of the edges of the cycle.) Finally, there are $5^{n-k}$ ways to choose the special points of vertices outside $C$, and then $(5n/2-k)!$ ways to pair up the points of appropriate types. Let  $a_{i}$ denote the number of orientations of $C$   with $i$ vertices of in-degree 2. Dividing the number of triples by $ M(5n) \E Y$,  and summing over all $i$, we obtain
$$
 \frac {\E  (Y X_k) }{\E  Y}\sim \sum_{0 \le i \le k/2} a_{i} \frac {  20^k [n]_k {n-2i \choose n/2-i}   3^{2i} 5^{n-k} (5n/2-k)!} {2k  {n \choose n/2} 5^n (5n/2)! } \sim \sum_{0 \le i \le k/2} \frac{a_{i}}{2k} \left( \frac 85 \right)^k \left( \frac 32 \right)^{2i}  .
$$
This gives~\eqn{mudef} with 
$$
\mu_k = \frac {1}{2k} \cdot \left( \frac 85 \right)^k \sum_{0 \le i \le k/2} a_{i}   \left( \frac 32 \right)^{2i}.
$$

 To find $a_{i}$, one can select the $2i$ vertices of $C$ that are to have  out-degree 0 or 2  in $C$, and after this there are exactly two ways to orient $C$. Hence $a_i=2{k\choose 2i}$, and this is the coefficient of $x^{2i}$ in 
 $
q(x):= 2 \left( 1 +   x \right)^k.
$ 
It follows that
$$
\sum_{0 \le i \le k/2} a_{i}  \left( \frac 32 \right)^{2i} = \frac 12 \Big( q(3/2)+q(-3/2) \Big) = \left( \frac 52 \right)^k + \left( - \frac 12 \right)^k,
$$
and thus
$$
\mu_k = \frac {1}{2k}    \big( 4^k  +  ( - 4/5)^k \big).
$$
 
As mentioned above, the derivation of~\eqn{mudef} is a straightforward generalisation of this. One starts with a set of cycles instead of one cycle, and following the same argument, the effects of those cycles lead to multiplicative factors in the counting that turn out to be asymptotically independent.  We omit the details.

The final step of small subgraph conditioning requires us to compute (see~\eqn{lambda})
$$
 \delta_k = \frac {\mu_k}{\lambda_k} - 1 = \left( - \frac 15 \right)^k 
$$
and then, using $-\log(1-x) = \sum_{k \ge 1} x^k / k$, 
$$
\exp \bigg( \sum_{k \ge 1} \lambda_k \delta_k^2 \bigg) = \exp \bigg( \frac 12 \sum_{k \ge 1} \frac {1}{k} \left( \frac {4}{25} \right)^k \bigg) = \exp \left( -\frac 12 \log \left( 1 - \frac {4}{25} \right) \right) = \frac {5}{\sqrt{21}}.
$$
The fact that this is coincides with the right hand side of~\eqn{ratio} implies, by~\cite[Theorem~4.1]{Wo}, that $\Prob(Y>0)\sim 1$. That is a random multigraph in $\mathcal{P}_{n,d}$ a.a.s.\ has a nowhere-zero flow over $\Z_3$. Theorem~\ref{thm:main}  now follows, in view of the comments at the start of this section.

\end{document}